\documentclass[11pt]{article}
\usepackage{amsmath,amsfonts,latexsym}
\usepackage{amscd,amssymb,amsmath}
\usepackage{graphicx}
\usepackage[english]{babel}
\usepackage{indentfirst,enumerate}
\usepackage[T1]{fontenc}
\usepackage[cp1250]{inputenc}
 \usepackage[colorlinks=true]{hyperref}
\hypersetup{urlcolor=blue, citecolor=red}

\newtheorem{theorem}{Theorem}

\newtheorem{proposition}{Proposition}

\newtheorem{definition}{Definition}

\newtheorem{remark}{Remark}
\newtheorem{lemma}{Lemma}
\newtheorem{corollary}{Corollary}

\newtheorem{conjecture}{Conjecture}

\newcommand{\hbx}{\hfill$\Box$}
\newcommand{\ep}{\varepsilon}

\begin{document}
\title{Firing map of an almost periodic input function}

\author{Wac{\l}aw Marzantowicz and Justyna Signerska}




\date{}
\maketitle


\medskip


\bigskip

\begin{abstract}
In mathematical biology and the theory of electric networks the firing map of an integrate-and-fire system is a notion of importance. In order to prove useful properties of this map authors of
 previous papers assumed that the stimulus function $f$ of the system $\dot{x}= f(t,x)$ is continuous and usually periodic in the time variable. In this work we show that the required properties of
 the firing map for the simplified model $\dot{x}=f(t)$ still hold if $f \in L_{loc}^1(R)$ and $f$ is an almost periodic function. Moreover, in this way we prepare a formal framework for next study of a
  discrete dynamics of the firing map arising from almost periodic stimulus that gives information on consecutive resets (spikes).
\end{abstract}
\section*{Introduction}
The integrate-and-fire models are used mainly in neuroscience to describe nerve-membrane voltage response to a given input. In these, usually one-dimensional models,
\begin{equation}\label{row ogolne}
\dot{x}=f(t,x), \;\;\; f:\mathbb{R}^2\to\mathbb{R}\,
\end{equation}
the continuous dynamics govern by a differential
equation is interrupted by threshold - reset behavior which is supposed to mimic spiking in real neurons. Precisely, the system evolves from an initial condition according to the differential equation as long as
the threshold-value of a dynamical variable is achieved. Then there is an immediate resetting to the resting-value and the dynamics continues from the new initial condition. The resting
and threshold values are typically set to $0$ and $1$, respectively. The sequence of consecutive spikes can be described by the iterates of the \emph{firing map}.

So far an analytical approach was carried mainly towards models with
periodic input, where the firing map is a lift of a degree one
circle map and the tools of rotation theory can be
used
(\cite{brette1}, \cite{Car-Ong}, \cite{gedeon}). Properties of the
firing map were also investigated with the help of numerical
simulations (\cite{mode},\cite{keener1}). In analytical survey it
was always assumed that the input function is smooth enough. In this
work we study the case when the function on the right hand side of a
differential equation is in general not continuous and not periodic
but almost periodic.
\section*{Properties of the firing map}\label{firing map}
\begin{definition} \rm
For the equation (\ref{row ogolne}) we define a map $\Phi: \mathbb{R} \to
\mathbb{R}$,
  $$\Phi(t)=\inf\{s>t: \
\, x(s;t,0)=1\}\,,$$
 where $x(\cdot;t,0)$ is a
solution of  (\ref{row ogolne}) satisfying the initial condition $(t,0)$.
\end{definition}
\noindent A natural domain of  $\Phi$ is the set $
D_{\Phi}=\{t\in \mathbb{R}: \quad \exists_{s>t} \ x(s;t,0)=1\}$. The mapping $\Phi$ assigns to each  $t\in D_{\Phi}$  the value $\Phi(t)$,
equal to the time after which a trajectory starting from $x=0$ at
the moment $t$ reaches the line  $x=1$ for the first time. This is the time
at which we have the \emph{firing} (alternatively called also a \emph{spike}). Immediately after that the solution
$x(t)$ is ''reset'' to zero and next the system given by
(\ref{row ogolne}) evolves from new starting point $(\Phi(t),0)$ until a
moment $\Phi^2(t)$ satisfying  $x(\Phi^2(t);\Phi(t),0)=1$. Thus the iterates $\Phi^n(t)$ of the firing map can be defined recursively as
\begin{displaymath}
\Phi^{n}(t)=\inf\{s>\Phi^{n-1}(t): \ x(s;\Phi^{n-1}(t),0)=1\}.
\end{displaymath}

The constant functions $x_r=0$ and $x_\tau=1$ are called  {\it the
reset} and  {\it  threshold} functions respectively. One can investigate varying threshold and reset functions (e.g.\cite{gedeon}) but the general cases reduces to the above (\cite{brette1}).

\subsection*{Perfect Integrator Model}

Consider the system called \emph{Perfect Integrator Model} (after \cite{brette1}),
i.e. the scalar ODE:
\begin{equation}\label{row1}
\dot{x}=f(t),
\end{equation}
with ''resetting mechanism''.
Throughout the rest of this paper we assume that $f\in L^1_{\textrm{loc}}(\mathbb{R})$, unless some additional assumptions are stated, and write simply $\int_{a}^{b}f(u)\,du$ meaning the Lebesque integral $\int_{A}f(u)\,du$ for $A=[a,b]$. We will denote $x(t;t_0,x_0):=\int_{t_0}^tf(u)\,du+x_0$ and call $x(\cdot;t_0,x_0)$ the solution of (\ref{row1})
(satisfying the initial condition $(t_0,x_0)$) or a trajectory of (\ref{row1}).
A locally integrable function has desired properties: the function $F(t)=\int_{t_0}^{t}f(u)\,du$ is a continuous function $F: (t_0,\infty)\to \mathbb{R}$ and it assures that the system (\ref{row1})
produces only a finite number of spikes in any finite interval, as we will see after Corollary \ref{uwagadopisana}. Moreover, this is a significantly larger class of functions than was studied
before (\cite{brette1}, \cite{Car-Ong}, \cite{mode}, \cite{gedeon}, \cite{keener1}): it includes, for instance, non-continuous piece-wise constant stimulus functions (e.g. the Haar wavelets), often used in
neuroscience (see e.g \cite{izykiewicz} and references therein).

  Since every two trajectories of (\ref{row1}) differ by a constant and the solution
$x(t;a,x_0)=\int_{a}^{t}f(u)\,du+x_0$ is continuous for  $f\in
L^{1}_{\textrm{loc}}(\mathbb{R})$,  we get the following
\begin{corollary}\label{uwagadopisana}
If $f\in L^{1}_{\textrm{loc}}(\mathbb{R})$, then for the equation
(\ref{row1}) the consecutive iterates of the firing map are equal to
\begin{equation}\label{szybkieiteracje}
\Phi^n(t)=\min\{s>t: \ x(s;t,0)=n\}
\end{equation}
and  there is only a finite number of firings in every bounded
interval.
\end{corollary}

Now we state necessary and sufficient conditions under which the
firing map $\Phi:\mathbb{R}\to\mathbb{R}$ is properly defined, it is
when $\textrm{D}_{\Phi}=\mathbb{R}$. In this case $\Phi$ induces a semi-dynamical system on
$\mathbb{R}$.
\begin{lemma}\label{lematdopisany}
For the model $\dot{x}=f(t)$
the firing map $\Phi:\mathbb{R}\to\mathbb{R}$ is well
defined if and only if
\begin{equation}\label{warunekkonidost}
\limsup_{t\to\infty}\int_{0}^{t}f(u) d u=\infty.
\end{equation}
In particular, for any $t_0\in\mathbb{R}$ (\ref{warunekkonidost}) is
the necessary and sufficient condition for all the iterations
$\{\Phi^n(t_0)\}_{n=1}^{\infty}$ to be finite real valued.
\end{lemma}
\textbf{Proof.} Choose $t_0\in\mathbb{R}$. We want all the iterations
$\{\Phi(t_0),\Phi^2(t_0),\Phi^3(t_0), ...\}$  to be well-defined
according to (\ref{szybkieiteracje}), that is we demand
that for all $n\geq 1$ there exists $s>\Phi^{n-1}(t_0)$ such that $x(s)=n$. However, this is equivalent to demanding that the solution $x:[t_0,\infty)\to\mathbb{R}$ with $x(t_0)=0$
is unbounded from above. But
\begin{displaymath}
\limsup_{t\to\infty}x(t)=\infty\iff\limsup_{t\to\infty}\int_{t_0}^{t}f(s)
d s=\infty\iff\limsup_{t\to\infty}\int_{0}^{t}f(s) d s=\infty.
\end{displaymath}
Obviously, $\limsup_{t\to\infty}\int_0^t f(s) d s=\infty$ asserts that
all the solutions of (\ref{row1}) are unbounded from above. \hbx
\vspace{0.5cm}

Throughout the rest of the paper we assume that (\ref{warunekkonidost}) is satisfied.
By using a property of the Lebesque integral which says that if the integral $\int_{A} f(t) d t$ of a non-negative function on a measurable set $A \subset \mathbb{R}$, $\mu(A)>0$, is positive then there exist $t_1 < t_2$ such that $\int_{t_1}^{t_2} f(t) dt >0$, one obtains
\begin{lemma}\label{monotonicity of solutions}
Any solution $x(t)$  of
the problem (\ref{row1}) is non-decreasing iff $f(t)\geq 0$ and is increasing iff $f(t)>0$ a.e. in $ \mathbb{R}$.
\end{lemma}

\begin{lemma}\label{monotonicity of firing} Let $f\in
L^{1}_{\textrm{loc}}(\mathbb{R})$ and $\Phi: \mathbb{R} \to \mathbb{R} $ be the firing map
associated with (\ref{row1}).

Then $\Phi$ is increasing,
  correspondingly,  non-decreasing, iff
$f(t)>0$, or $f(t)\geq 0$ respectively, a.e. in $\mathbb{R}$.
\end{lemma}
\textbf{Proof.} The proof that $f(t)\geq 0$ ($f(t)>0$) a.e. induces non-decreasing (increasing) firing map $\Phi$ is almost immediate from the previous lemma. Conversely, suppose that there
 exists a non-zero measure set $A$ such that $f(t)<0$ for $t\in A$. From Lemma \ref{monotonicity of solutions} there exist $t_1,t_2$, $t_1<t_2$, such
that $x(t_2)-x(t_1)=\int_{t_1}^{t_2}f(u)\, du<0$ where $x(\cdot)$ is an arbitrary solution of (\ref{row1}). Suppose that $\Phi(t_1)\geq t_2$. Then $1=\int_{t_1}^{t_2}f(u)\, du+\int_{t_2}^{\Phi(t_1)}f(u)\, du$
implies $\int_{t_2}^{\Phi(t_1)}f(u)\, du>1.$ Consequently $\Phi(t_2)<\Phi(t_1)$. If $\Phi(t_1)<t_2$, let then $n\in\mathbb{N}$ be such that $\Phi^{n}(t_1)<t_2$ but $\Phi^{n+1}\geq t_2$ (existence of such $n$ is ensured by the fact that there can be only a finite number of firings in the interval $[t_1,t_2]$). Then, as before, we obtain that $\Phi(t_2)<\Phi(t_*)$, where $t_*=\Phi^{n}(t)<t_2$ because in this case $\int_{t_*}^{t_2}f(u)\,du<0$. Analogously we show that when $\Phi$ is strictly increasing, then $f$ is positive, perhaps with the exception of some zero-measure set.
\hbx
\vspace{0.5cm}

Some general results concerning continuity and injectivity of the
firing map for the models of the type $\dot{x}=f(t,x)$ when $f$ is
smooth enough are given in \cite{Car-Ong}. However, for the Perfect
Integrator we obtain a detailed description of continuity if only $f\in L^{1}_{\textrm{loc}}(\mathbb{R})$.

Let $f\in L^{1}_{\textrm{loc}}(\mathbb{R})$, $f(t)\geq 0$ almost everywhere. For $a\in \mathbb{R}$
consider $\Phi^{-1}(a) \subset \mathbb{R}$. Let $\bar{a}= \sup\{ t\in
\Phi^{-1}(a)\}$. Note that $\bar{a} < \infty$, because $\bar{a}<a$ since for any $t\in\mathbb{R}$ $t<\Phi(t)$. Moreover, there exists an
interval $[\bar{a}, \bar{a} +\bar{\delta}]$ such
 that for every $t$ satisfying $\bar{a} < t \leq \bar{\delta}$ we have
 $\int_{\bar{a}}^t f(u) \,du >0$, or equivalently $\mu(\{u: \, f(u)>0\} \cap [\bar{a}, \bar{a} + t]) >0$.
 Indeed, if there exists $t$ such that $a> t>\bar{a}$ and $\int_{\bar{a}}^tf(u)\, du = 0$
 then $ \int_{t}^{a} f(u)\, du  = \int_{\bar{a}}^{a} f(u)\,du
 =1$ and hence $t\in\Phi^{-1}(a)$ contrary to the definition
 of $\bar{a}$.

\begin{proposition}\label{discontunity}
If $f\in L^{1}_{\textrm{loc}}(\mathbb{R})$ and $f(t)\geq 0$ almost
everywhere. Then
\begin{itemize}
\item[i)] {$\Phi$ is left continuous,}
\item[ii)] {$\Phi$  is not right continuous at every point $\bar{a} \in \Phi^{-1}(a)$ for which there exists $\delta_0>0$ such that $f(t)=0$ almost everywhere in $[a,
a+\delta_0]$. Moreover such points are the only points of
discontinuity of $\Phi$. }
\end{itemize}
\end{proposition}
\textbf{Proof.} For the first part of the statement take a
sequence $t_n \to t_0$, $t_n <t_0$. Then $1 =\int_{t_n}^{\Phi(t_n)}
f(u)\, du $ and $1 = \int_{t_0}^{\Phi(t_0)} f(u)\,du $ implies
$\int_{t_n}^{t_0} f(u) \,du = \int_{\Phi(t_n)}^{\Phi(t_0)} f(u)\, du
$. It follows that $\int_{\Phi(t_n)}^{\Phi(t_0)} f(u)\, du \to
0$ if $n\to \infty$. We have to prove that $\Phi(t_n) \to
\Phi(t_0)$. Note that since $\Phi$ is non-decreasing, $\Phi(t_n)\leq\Phi(t_0)$. If for some $\delta_0>0$ $f(t) = 0 $ almost everywhere in
the interval $[\Phi(t_0) - \delta_0,\Phi(t_0)]$ then $
\int_{t_0}^{\Phi(t_0) - \delta_0} f(u)\, du=1$ which gives a contradiction to the definition of $\Phi$. Thus
$f(t)>0 $ on a  set of  positive  measure  in any small interval
$[\Phi(t_0) - \delta,\Phi(t_0)]$ and consequently $\Phi(t_n)
\to \Phi(t_0)$.

To prove the second part of the statement let us take $t_n\to \bar{a},
\; t_n>\bar{a}$. We have already showed that $\int_{\bar{a}}^{t_n}
f(u)
 \,du \, >0 $. Then
$\Phi(t_n) \geq \Phi(\bar{a}) +\delta_0 $ for large $n$, because for
every $b$, $ \Phi(\bar{a}) \leq b \leq \Phi(\bar{a}) +\delta_0$ we have $$
\int_{t_n}^b f(u)\,du =  \int_{t_n}^{\Phi(\bar{a})} f(u)
 \,du =\, \int_{\bar{a}}^{\Phi(\bar{a})} f(u)\,du - \int_{\bar{a}}^{t_n} f(u)
 \,du \, < 1.$$ It follows that $\Phi(\bar{a}_+) \geq \Phi(\bar{a}) +\delta_0$ and consequently $\Phi$ is not right  continuous  at $\bar{a}$.

Conversely, let $t_n \to t_0$,  $t_n> t_0$ and $\Phi(t_n) \geq
\Phi(t_0) +\bar{\delta}$. Then $\int_{t_0}^{t_n} f(u)\, du >0$,
since otherwise $1= \int_{t_0}^{\Phi(t_0)} f(u)\, du =
\int_{t_n}^{\Phi(t_0)} f(u)\, du $, which gives
$\Phi(t_n)=\Phi(t_0)$ contrary to the supposition. Next, if for
every $\delta_0 >0$ we have  $\mu(\{f(u)>0\}
\cap [\Phi(t_0), \Phi(t_0) +\delta_0]) >0$, then
$\int_{\Phi(t_0)}^{\Phi(t_0)+\delta_0} f(u)\, du >0$. We can take
$\delta_0 < \bar{\delta}$. For large $n$ $\int_{t_0}^{t_n}f(u)\, du $ is
small, say smaller than $\int_{\Phi(t_0)}^{\Phi(t_0)+\delta_0} f(u)\,
du$. This gives $ \int_{t_n}^{\Phi(t_0) +\delta_0} f(u)\,du$ $=
\int_{t_0}^{\Phi(t_0)} f(u) \,du -\int_{t_o}^{t_n} f(u)\,du +
\int_{\Phi(t_0)}^{\Phi(t_0) +\delta_0} f(u)\,du >1$, thus $\Phi(t_n) <
\Phi(t_0) + \delta_0$ contrary to our supposition that $\Phi(t_n) \geq
\Phi(t_0) +\bar{\delta}$.
\hbx
\vspace{0.5cm}

As a direct consequence of Proposition \ref{discontunity} we get the
following:

\begin{corollary}\label{continuity}
If $f\in L^{1}_{\textrm{loc}}(\mathbb{R})$
and $f(t)>0$ almost everywhere, then $\Phi$ is continuous.
\end{corollary}

\section*{Firing rate}

\begin{definition}\label{def_firing}
For every $t\in \mathbb{R}$ one defines  the \emph{firing rate}, denoted by
$r(t)$, as
\begin{equation}\label{row2}
r(t):=\lim_{n\to\infty}\frac{n}{\Phi^{n}(t)}.
\end{equation}
\end{definition}

 In general for an equation  $\dot{x}=f(t,x)$ the limit of
(\ref{row2}) might not exist for some $t$. Moreover, even if the
limit exists for all $t$, it might depend on $t$. However, for the simplified model (\ref{row1}) we have the
following theorem, which was proved in \cite{brette1} and the proof is valid for $f\in L^{1}_{\textrm{loc}}(\mathbb{R})$:
\begin{theorem}\label{tw1}
Suppose that for the model  (\ref{row1}) there exists a finite
limit
\begin{equation}\label{row3}
r=\lim_{t\to\infty}\frac{1}{t}\int_{0}^{t}f(u)d u.
\end{equation}
Then for every point $t_0\in\mathbb{R}$ the firing rate $r(t_0)$
exists and is given by the formula (\ref{row3}). In particular the
firing rate $r(t)$ does not depend on $t$.
\end{theorem}
\subsection*{Almost periodic functions}
 In next we consider the model (\ref{row1}) with the
 function
 $f:\mathbb{R}\to\mathbb{R}$ almost periodic. There are various classes of almost periodic functions, e.g. in the sense of Bohr, Stepanov, Weyl and Besicovitch. They are defined as
 the closure of the space of generalized trigonometric polynomials under corresponding norms. The notion of almost periodicity generally holds for functions $f:\mathbb{R}\to\mathbb{C}$.
 However, we will only consider almost periodic functions taking real values.
\begin{definition}\label{Bohr alm per}
A continuous function $f:\mathbb{R}\to\mathbb{C}$ is \emph{almost periodic function in the sense of Bohr} (alternatively called also uniformly almost periodic) if for any $\ep>0$ the set $E\{\ep,f(t)\}$ of
all the numbers $\tau$ such that for all $t\in\mathbb{R}$ $|f(t+\tau)-f(t)|<\ep$ is relatively dense, i.e. there exists $l_{\ep}>0$ such that in any interval of length $l_{\ep}$ there is at least one $\tau$
satisfying this inequality.

The set $E\{\ep,f(t)\}$ is called the set of $\ep$-almost periods of $f$.
\end{definition}
Equivalently we can say that $f$ is Bohr almost periodic if it is the limit of a sequence of generalized trigonometric polynomials $P_n(t)=\sum_{i=1}^{k(n)}c_i\exp(\imath \lambda_i t)$, $c_i\in\mathbb{C}$ and $\lambda_i\in\mathbb{R}$,
under the norm $\|f\|=\sup_{t\in\mathbb{R}}|f(t)|$. Any Bohr almost periodic function is uniformly continuous and bounded (\cite{bezykowicz},\cite{corduneanu}).
\begin{definition}\label{stepanov alm per}
A function $f:\mathbb{R}\to\mathbb{C}$, $f\in L^{p}_{\textrm{loc}}(\mathbb{R})$, is \emph{Stepanov almost periodic} if for any $\ep>0$ the set $SE\{\ep,f(t)\}$ of all the numbers $\tau$ such
that $\|f(t+\tau)-f(t)\|_{\textrm{St.},r,p}<\ep$ is relatively dense, where
$$\|f\|_{\textrm{St.},r,p}:=\sup_{t\in\mathbb{R}}[\frac{1}{r}\int_{t}^{t+r}|f(u)|^p\,du]^{1/p}, \quad r>0, \ 1\leq p<\infty.$$
The set $SE\{\ep,f(t)\}$ is called the set of $\ep$-Stepanov almost periods of $f$.
\end{definition}
The space of Stepanov almost periodic functions can be obtained also as the closure of generalized trigonometric polynomials under the norm $\|\cdot\|_{\textrm{St.},r,p}$. Any Stepanov almost periodic function
which is uniformly continuous is Bohr almost periodic. Different values of $r$ in $\|\cdot\|_{\textrm{St.},r,p}$ give obviously different norms but the same topology.
Thus we will consider Stepanov almost periodic functions with $r=1$. If $f$ is a $p_1$-Stepanov almost periodic function, then it also $p_2$-Stepanov almost periodic for $p_2<p_1$ with fixed $r$ (\cite{andres},
 \cite{bezykowicz}). One can also consider almost periodic functions in the sense of Weyl or Besicovitch, where the class of Besicovitch almost periodic functions with $p=1$ is larger and includes all
 the other mentioned types of almost periodicity. For an almost periodic function $f:\mathbb{R}\to\mathbb{C}$ of any of these types there exists a finite \emph{mean} (\cite{andres},\cite{bezykowicz}, \cite{corduneanu}):
\[
\mathcal{M}\{f(t)\}=\lim_{t\to\infty}\frac{1}{t}\int_{0}^{t}f(u)\,du.
\]
From this we conclude that:
\begin{corollary}\label{maly wniosek1}
For the model (\ref{row1}), where $f:\mathbb{R} \to \mathbb{R}$ is almost periodic in any of the above sense the firing rate exists and equals the mean of a function, provided that the mean is nonnegative
(otherwise the firing rate is zero).
\end{corollary}
Note that if $\mathcal{M}\{f(t)\}<0$, then (\ref{warunekkonidost}) does not hold and there is no firing.
\begin{definition}\label{def limit periodic}
A function $f:\mathbb{R}\to\mathbb{C}$ is limit-periodic (in the sense of Bohr) if it is a limit of uniformly converging sequence $\{f_n\}$ of continuous periodic functions.
\end{definition}
Note that generalized trigonometric polynomials are in general not periodic (consider, for example, $P(t)=\sin(\sqrt{2}t)+\sin(2t)$) and that a limit-periodic function is an uniformly almost periodic function. One can read in \cite{andres} in details about the theory of almost periodic functions.
 We will write shortly that $f$ is u.a.p. or S.a.p., meaning uniform (Bohr) or Stepanov (with $r=1$ and $p=1$) almost periodicity, respectively.

 The mean $\mathcal{M}$ is continuous \cite{bezykowicz}, meaning that $\mathcal{M}\{f_n\}\to\mathcal{M}\{f\}$ if $f_n\to f$ in any of the ``almost-periodic norms" (i.e. the norms with respect to which a given space
 of almost periodic functions can be obtained as the closure of the space of generalized trigonometric polynomials). Thus the following holds:
\begin{theorem}\label{firing rate for almost periodic}
Let $(f_n)_{n=1}^{\infty}$ be a sequence of
almost-periodic functions $f_n:\mathbb{R}\to\mathbb{R}$ of the Bohr, Stepanov, Weyl or Besicovitch type converging to $f$
in uniform, Stepanov, Weyl or Besicovitch norm, respectively. Let  $r_n$ be the firing rate of
the equation $\dot{x}=f_n(t)$,  for $n=1,2,3,...$, and $r$ the firing rate for (\ref{row1}).

Then $r=\lim_{n\to\infty}r_n$.
\end{theorem}

We end this section with propositions showing when the
existence of the positive mean value of $f$ is a necessary and sufficient condition for existence of the
firing map $\Phi:\mathbb{R}\to\mathbb{R}$. By elementary calculations one proves
\begin{proposition}\label{mean implies firing}
Let $f\in L^{1}_\textrm{loc}(\mathbb{R})$. If
$\mathcal{M}\{f\}>0$ then
$ \limsup \int_0^{t} f(u)\,du =
\infty$. Consequently the firing map $\Phi$ of (\ref{row1}) is well-defined  for
all $t\in \mathbb{R}$.

Conversely, if $f(t)= \sum _{j=1}^k  c_j\exp(\imath \lambda_j \, t)$, then $\limsup_{t\to\infty} \int_0^t
f(u)\, du =\infty $ implies $ \mathcal{M}\{f\}>0$.
\end{proposition}
In the second case it must hold that $\mathcal{M}\{f\}=c_{j_*}>0$, where $c_{j_*}$ corresponds to $\lambda_{j_*}=0$.

The mean has the following property, which in general does not hold for functions that have finite $\sup$- or Stepanov norms but are not almost periodic:
\begin{lemma}\label{o sredniej dla nieujemnej}
Let $f:\mathbb{R}\to\mathbb{R}$ be a real non-negative (a.e.) u.a.p. (/S.a.p.) function. If $\mathcal{M}\{f(t)\}=0$, then $f(t)=0$ for  all (/almost all) $t\in\mathbb{R}$.
\end{lemma}
\textbf{Proof.} The proof of the statement for u.a.p. functions can be found in \cite{bezykowicz}. We will prove the statement for $f$ being a S.a.p. function. Suppose
on the contrary that there exists a non-zero measure set $A$ such that for all $t\in A$ we have $f(t)>0$. Then there exists $a\in\mathbb{R}$ such that $
\int_{a}^{a+1}f(u)\,du>\lambda$
for some $\lambda>0$ as follows from Lemma \ref{monotonicity of firing}. Since $f$ is S.a.p. the set of all the $\tau$ such that
for all $t\in\mathbb{R}$ $\int_{t}^{t+1}|f(u+\tau)-f(u)|<\frac{\lambda}{2}$ is relatively dense.
It follows that there exists $l_{\frac{\lambda}{2}}>0$ such that in any interval $I$ of length $l_{\frac{\lambda}{2}}$ we can find $\tau\in I$ such that
$\int_{a+\tau}^{a+\tau+1}f(u)\,du=\int_{a}^{a+1}f(u+\tau)\,du>\frac{\lambda}{2}$: Indeed, if $\int_{a}^{a+1}|f(u+\tau)-f(u)|\,du<\frac{\lambda}{2}$ and $\int_{a}^{a+1}f(u)\,du>\lambda$, then $\frac{\lambda}{2}<\int_{a}^{a+1}f(u)\,du-\int_{a}^{a+1}|f(u+\tau)-f(u)|\,du\leq\int_{a}^{a+1}f(u+\tau)\,du$.
We can assume that $l_{\frac{\lambda}{2}}>1$. Then we can find a sequence $\{\tau_n\}$ where $\tau_{n}\in[2n l_{\frac{\lambda}{2}},2(n+1)l_{\frac{\lambda}{2}}-1]$ such that
\[
\int_{a+2kl_{\frac{\lambda}{2}}}^{a+2(k+1)l_{\frac{\lambda}{2}}}f(u)\,du>\int_{a+\tau_{k}}^{a+\tau_{k}+1}f(u)\,du>\frac{\lambda}{2}
\]
for any $k\in\mathbb{N}\cup\{0\}$ because $a+\tau_{k},a+\tau_{k}+1\in[a+2kl_{\frac{\lambda}{2}},a+2(k+1)l_{\frac{\lambda}{2}}]$. Consequently,
\begin{equation}
\begin{split}
0=\mathcal{M}\{f\}
=&\lim_{t\to\infty}\frac{1}{t}\int_{0}^{t}f(u)\, du=\lim_{t\to\infty}\frac{1}{t}\int_{a}^{a+t}f(u)\,
 du\nonumber\\
 =&\lim_{n\to\infty}\frac{1}{2nl_{\frac{\lambda}{2}}}\int_{a}^{a+2nl_{\frac{\lambda}{2}}}f(u)\,du>\frac{1}{2nl_{\frac{\lambda}{2}}}\,n\,\frac{\lambda}{2}\,=\,\frac{\lambda}{4l_{\frac{\lambda}{2}}}
\nonumber
\end{split}
\end{equation}
and we obtain a contradiction.
\hbx
\vspace{0.5cm}

\begin{proposition}\label{mean and firing dla nonnegative}
Suppose that $f(t)\geq 0$ (a.e.) is an u.a.p. (/S.a.p.) function. Then $\mathcal{M}\{f\}>0$ if and only if the firing map is well-defined.
\end{proposition}
\textbf{Proof.} Sufficiency was covered in the previous proposition, necessity follows from Lemma \ref{o sredniej dla nieujemnej}. Indeed, if $\mathcal{M}\{f\}=0$ for non-negative (a.e.) u.a.p. (/S.a.p.) function $f$, then $f(t)=0$
for all $t\in\mathbb{R}$ (/a.e. in $\mathbb{R}$). In both cases the condition $\limsup_{t\to\infty}\int_{0}^t f(u)\,du=\lim_{t\to\infty}\int_{0}^t f(u)\,du=\infty$ is not satisfied.
\hbx
\vspace{0.5cm}

\begin{remark}\label{welvets}
The assumption that  a stimulus function  $f\in
L^{1}_{\textrm{loc}}(\mathbb{R})$ allows us also to take as $f$ a
linear combination of the Haar wavelets which is more natural for many
problems with discontinuous stimulus function.
\end{remark}
\section*{Displacement map}
\begin{definition}\label{displacement map}{\rm
Let $\Phi: \mathbb{R} \to \mathbb{R}$ be a map of the real line. The {\em displacement map} $\Psi: \mathbb{R} \to \mathbb{R} $ of $\Phi$ is defined
as
$$ \Psi(t)\,:= \Phi(t)-t \,.$$}
\end{definition}

If $\Phi$ is a firing map then the displacement map $\Psi(t)$ says how long we have to wait for a next firing if we know that there was a firing at the time $t$.

Firstly, we will consider the firing map and its displacement for the general case, i.e. for the equation (\ref{row ogolne}). The following observation was first made in \cite{keener1}:
\begin{proposition}\label{displacement okresowe}
If the function $f$ in (\ref{row ogolne}) is periodic in $t$ (that is,
there exists $T$ such that for all $x$ and $t$ we have
$f(t,x)=f(t+T,x)$), then the firing map $\Phi$ has periodic
displacement. In particular for $T=1$ we have $\Phi(t+1)=\Phi(t)+1$ and thus $\Phi$ is a lift of a degree one circle map under the standard projection $\mathfrak{p}: t\mapsto \exp(2\pi\imath t)$.
\end{proposition}

Therefore for the models (\ref{row ogolne}) with $f$ periodic in $t$, the firing map $\Phi$ has periodic displacement and induces a degree one circle map (this map is sometimes in the literature
referred to as the \emph{firing phase map}). In this case the tools of the rotation theory (\cite{misiurewicz}) can be used to study the properties of $\Phi$. Especially, the rotation number
(which is simply the reciprocal of the firing rate, if the unique non-zero firing rate exists) or the rotation intervals are closely related to the so-called phase-locking phenomena (\cite{brette1}, \cite{keener1}).

There arises a natural question of whether for the map $f$ almost periodic in $t$ we obtain the firing map with almost periodic displacement. We will tackle this problem for the Perfect Integrator Model.

 Firstly, let us consider the case of limit-periodic input:
\begin{theorem} \label{twierdzenie o limitperiodic} Let $\Phi:\mathbb{R}\to \mathbb{R}$ be the firing map
induced by (\ref{row1}). If $f:\mathbb{R}\to\mathbb{R}$ is limit periodic and $f(t)>\delta$ for some $\delta>0$, then the firing map
$\Phi$ has limit periodic displacement.
\end{theorem}
\textbf{Proof.} Denote the displacement of $\Phi$ by $\Psi$. Suppose, without the loss of generality, that $0<\delta<1$ and let $\ep>0$ be any small number. We can assume that
$\delta>\ep/2>0$. Since $f$ is limit periodic, there exists
continuous periodic function $\widetilde{f}:\mathbb{R}\to\mathbb{R}$ such that $ |\widetilde{f}(t)-f(t)|<\delta^2\ep/4$ for all $t$.
Choose the initial condition $(t_0,0)$. By $x(t)$ denote
the solution of (\ref{row1}) satisfying
$x(t_0)=0$ and by $\widetilde{x}(t)$ the solution of
\begin{equation}\label{eq5}
\dot{x}=\widetilde{f}(t)
\end{equation}
with the same initial condition.
Equation (\ref{eq5}) gives rise to the firing map $\widetilde{\Phi}(t)$ with the displacement $\widetilde{\Psi}(t)$. We will show that $|\Psi(t_0)-\widetilde{\Psi}(t_0)|<\ep$, which is equivalent
to $|\Phi(t_0)-\widetilde{\Phi}(t_0)|<\ep$, and then the theorem will follow from the fact that the firing map $\widetilde{\Phi}$ of periodically forced
Perfect Integrator Model has periodic displacement, which is
covered by Proposition \ref{displacement okresowe}. Since $\widetilde{f}$ approximates $f$ in the
uniform norm, we have that for all $t$ $\widetilde{f}(t) \geq \delta/2$.
We have to consider two possibilities: $\Phi(t_0)>\widetilde{\Phi}(t_0)$ and $\widetilde{\Phi}(t_0)>\Phi(t_0)$. Suppose that $\Phi(t_0)>\widetilde{\Phi}(t_0)$, i.e. the solution $\widetilde{x}(t)$
fires before $x(t)$. Then $x(\Phi(t_0))=1$ and $x({\widetilde{\Phi}}(t_0))<1$ by the definition of $\widetilde{\Phi}$ and $\Phi$. From the Mean Value
Theorem
\begin{equation}\label{eq7}
1-x(\widetilde{\Phi}(t_0))=f(\alpha)(\Phi(t_0)-\widetilde{\Phi}(t_0))
\end{equation}
for some $\alpha\in(\widetilde{\Phi}(t_0),\Phi(t_0))$. Now define $y(t):=\widetilde{x}(t)-x(t)$. Then $y(t)-y(t_0)=y^{\prime}(\xi_t)(t-t_0)$
for any $t>t_0$ and some $\xi_t\in(t_0,t)$. Since $y(t_0)=0$ and $y^{\prime}(\xi_t)=\widetilde{x}^{\prime}(\xi_t)-x^{\prime}(\xi_t)=\widetilde{f}(\xi_t)-f(\xi_t)\in(-\delta^2\ep/4, \delta^2\ep/4)$, we get
\begin{equation}\label{eq8}
\widetilde{x}(\widetilde{\Phi}(t_0))-x(\widetilde{\Phi}(t_0))<(\widetilde{\Phi}(t_0)-t_0)\frac{\delta^2\ep}{4}.
\end{equation}
Equations (\ref{eq7}) and (\ref{eq8}) yield $\Phi(t_0)-\widetilde{\Phi}(t_0)<\frac{ \delta^2\ep}{4 f(\alpha)}(\widetilde{\Phi}(t_0)-t_0)$, where the difference  $\widetilde{\Phi}(t_0)-t_0\leq\frac{2}{\delta}$
 because $1=\widetilde{x}(\widetilde{\Phi}(t_0))-\widetilde{x}(t_0)=\widetilde{f}(\beta)(\widetilde{\Phi}(t_0)-t_0)$ for some $\beta\in(t_0,\widetilde{\Phi}(t_0))$. Finally,
$\Phi(t_0)-\widetilde{\Phi}(t_0)<\frac{
\delta^2\ep}{4 f(\alpha)\widetilde{f}(\beta)}< \ep$ since $f(\alpha),\widetilde{f}(\beta)>\frac{\delta}{2}$. Analogously we prove that $\widetilde{\Phi}(t_0)-\Phi(t_0)<\ep$ for $\widetilde{\Phi}(t_0)>\Phi(t_0)$.
\hbx
\vspace{0.5cm}

From the proof of Theorem \ref{twierdzenie o limitperiodic} follows more general
\begin{proposition}\label{twierdzenie o przyblizaniu}
If $f:\mathbb{R}\to\mathbb{R}$, $f(t)>\delta>0$ is a limit of uniformly convergent sequence of continuous functions $f_n:\mathbb{R}\to\mathbb{R}$ inducing firing maps $\Phi_n$ with displacements $\Psi_n$,
the firing map $\Phi$ and the displacement $\Psi$ of $f$ can be obtained as the uniform limits of $\Phi_n$ and $\Psi_n$, respectively.
\end{proposition}
The above means that if the system $\dot{x}=f_{\gamma}(t)$ depends on the parameter $\gamma\in\mathbb{R}$ continuously with respect to the uniform norm, i.e. $\sup_{t\in\mathbb{R}}|f_{\gamma_1}(t)-f_{\gamma_2}(t)|$
is arbitrary small provided that the difference $|\gamma_1-\gamma_2|$ is small enough, then the corresponding firing maps $\Phi_{\gamma}(t)$ also change continuously with respect to the uniform topology. Furthermore we show that

\begin{theorem}\label{new almost periodic displacement}
Let $f:\mathbb{R}\to\mathbb{R}$ be a S.a.p. function such that $f(t)>\delta$ almost everywhere for some $\delta>0$. Then the firing map $\Phi$ induced by the equation (\ref{row1}) has u.a.p. displacement.
\end{theorem}
 \textbf{Proof.} We are to show that for any $\ep>0$ the set $E\{\ep,\Psi(t)\}$ of all the numbers $\tau\in\mathbb{R}$ such that $\sup_{t\in\mathbb{R}}|\Psi(t+\tau)-\Psi(t)|=\sup_{t\in\mathbb{R}}|\Phi(t+\tau)-\tau -\Phi(t)|<\ep$ is relatively dense.
Without the loss of generality we assume that $\delta<1$. Choose then $\ep>0$ and let $\tau\in SE\{\frac{\delta^2\ep}{2},f(t)\}$. Take $t\in\mathbb{R}$. Then by the assumption on $f$ $\int_{t}^{t+1}|f(u+\tau)-f(u)|\,du<\frac{\delta^2 \ep}{2}$.
Suppose that $\min\{\Phi(t),\Phi(t+\tau)-\tau\}=\Phi(t)$. By the definition of the firing map $\int_{t+\tau}^{\Phi(t+\tau)}f(u)\,du=\int_{t}^{\Phi(t)}f(u)\,du$ and it follows
that
\begin{equation}
\begin{split}
  0 =& \int_{t+\tau}^{\Phi(t+\tau)}f(u)\,du-\int_{t}^{\Phi(t)}f(u)\,du=
  \int_{t}^{\Phi(t+\tau)-\tau}f(u+\tau)\,du-\int_{t}^{\Phi(t)}f(u)\,du=\nonumber
   \\
   =& \int_{t}^{\Phi(t)}f(u+\tau)-f(u)\,du+\int_{\Phi(t)}^{\Phi(t+\tau)-\tau}f(u+\tau)\,du.\nonumber
  \end{split}
\end{equation}
Thus $|\int_t^{\Phi(t)}f(u+\tau)-f(u)\,du|=|\int_{\Phi(t)}^{\Phi(t+\tau)-\tau}f(u+\tau)\,du|$.
Since \nolinebreak[10]$\tau\nolinebreak[10]\in\nolinebreak[10]SE\{\frac{\delta^2\varepsilon}{2}, f(t)\}$,
$
  |\int_{t}^{\Phi(t)}f(u+\tau)-f(u)\,du| \leq \int_{t}^{\Phi(t)}|f(u+\tau)-f(u)|\,du\leq \int_{t}^{t+k}|f(u+\tau)-f(u)|\,du < k  \delta^2 \ep/2,
$
where $k\in\mathbb{N}$ is the smallest integer such that $\Phi(t)\leq t+k$. However, as $f(u)>\delta$ almost everywhere, $\int_{\Phi(t)}^{\Phi(t+\tau)-\tau}f(u+\tau)\,du>\delta(\Phi(t+\tau)-\tau-\Phi(t))$. Finally
 \[
 \Phi(t+\tau)-\tau-\Phi(t)<\frac{k \delta\ep}{2}<\frac{(\frac{1}{\delta}+1)\delta\ep}{2}=\frac{(1+\delta)\ep}{2}<\ep
 \]
because $\Phi(u)-u<1/\delta$ for any $u\in\mathbb{R}$, which can be estimated as in the proof of Theorem \ref{twierdzenie o limitperiodic}, and thus $k<(1/\delta+1)$. If $\min\{\Phi(t),\Phi(t+\tau)-\tau\}=\Phi(t+\tau)-\tau$ in a similar way we obtain that
$\Phi(t)-\Phi(t+\tau)+\tau<k \delta\ep/2$ where $k$ is the smallest integer such that $\Phi(t+\tau)-(t+\tau)<k$. Thus in any case $|\Phi(t+\tau)-\Phi(t)-\tau|<\ep$. It follows that $\tau\in E\{\ep,\Psi(t)\}$. Consequently the set $E\{\ep,\Psi(t)\}$ is relatively dense because it contains relatively dense set $SE\{\frac{\delta^2\varepsilon}{2},f(t)\}$.
\hbx
\vspace{0.5cm}

Note that if the displacement $\Psi(t)=\Phi(t)-t$ is uniformly almost periodic then it is uniformly continuous and the following conclusion is immediate:
\begin{corollary}
Under the assumptions of Theorem \ref{new almost periodic displacement} the firing map $\Phi:\mathbb{R}\to\mathbb{R}$ is uniformly continuous.
\end{corollary}
\begin{remark}\label{Kwapiszwork}
Rotation numbers and sets for maps of the real line
with almost periodic displacement were investigated by J. Kwapisz in
\cite{kwapisz}. Observe that Theorems \ref{twierdzenie o
limitperiodic} and \ref{new almost periodic displacement} ensure that the displacement of a firing
map is almost periodic if $f$ is so. On the other hand if a given almost periodic function is a firing map then the thesis of main theorems
 of \cite{kwapisz} can be deduced easier by a direct
argument. Anyway, it would give an opportunity to compare
geometrical definitions of rotation intervals of \cite{misiurewicz},
used in \cite{kwapisz}, with analytic formulas derived here.
\end{remark}

\section*{Prospects}

As we have already  said, our aim was  also to establish a framework
for study dynamical properties of the firing map $\Phi$ and its
displacement $\Psi$.

Now we state a conjecture that the differences of
consecutive spikes form an asymptotically semi-periodic, or respectively almost
periodic sequence if the input function is so (cf. \cite{berg} for definitions of
semi-periodic, and almost-periodic sequence).

Let $\Phi: \mathbb{R} \to \mathbb{R}$ be the firing map of the
equation $\dot{x}=f(t)$ with $f\in L^{1}_{\textrm{loc}}(\mathbb{R})$
and $\Psi:\mathbb{R}\to \mathbb{R}$ its displacement.
 For any $t\in \mathbb{R}$,  we consider a sequence
 $$ \eta_n(t)= \Phi^n(t) - \Phi^{n-1}(t) = \Psi(\Phi^{n-1}(t))
 $$
\begin{conjecture}
If $f$ is periodic, correspondingly almost-periodic, then for every
$t$ the sequence $\eta_n(t)$ is asymptotically semi-periodic,
respectively asymptotically almost-periodic.
\end{conjecture}


\medskip
FACULTY OF MATHEMATICS AND COMPUTER SCI.,\, ADAM MICKIEWICZ
UNIVERSITY OF POZNA{\'N}, ul. Umultowska 87,\, 61-614 Pozna{\'n},
Poland\\
\emph{E-mail address}: marzan@amu.edu.pl\\

\vspace{0.4cm}
\noindent INSTITUTE OF MATHEMATICS, POLISH ACADEMY OF SCIENCES, ul.\'Sniadeckich 8, 00-956 Warszawa, Poland\\
\emph{E-mail address}: j.signerska@impan.pl

\end{document}